\documentclass[a4paper,11pt]{amsart}

\usepackage{graphicx}
\usepackage{mathptmx}
\usepackage{amsmath}
\usepackage{amssymb}
\usepackage{enumitem}
\usepackage{xcolor}
\usepackage{epsfig}

\usepackage[FIGTOPCAP]{subfigure}
\usepackage{subfigure}

\DeclareMathOperator{\sech}{sech}

\newmuskip\pFqmuskip

\newcommand*\pFq[6][8]{%
  \begingroup 
  \pFqmuskip=#1mu\relax
  \mathcode`=\string"8000
  \begingroup\lccode`\~=`\,
  \lowercase{\endgroup\let~}\pFqcomma
  F^{#2}_{#3}{\left(\genfrac..{0pt}{}{#4}{#5}\bigg|#6\right)}%
  \endgroup
}
\newcommand{\pFqcomma}{\mskip\pFqmuskip}

\newtheorem{theorem}{Theorem}[section]
\newtheorem{lemma}[theorem]{Lemma}
\newtheorem{corollary}[theorem]{Corollary}

\begin{document}

\title[]{Probabilistic proof of a summation formula}

\author{Taekyun  Kim}
\address{Department of Mathematics, Kwangwoon University, Seoul 139-701, Republic of Korea}
\email{tkkim@kw.ac.kr}
\author{Dae San  Kim}
\address{Department of Mathematics, Sogang University, Seoul 121-742, Republic of Korea}
\email{dskim@sogang.ac.kr}



\begin{abstract}
The aim of this paper is to derive a summation formula for the series $\sum_{k=0}^{\infty} \frac{(-1)^{k}} {(2k+1)^{2n+1}}$ and an expression for $\zeta(2n+2)$ by using hyperbolic secant random variables.
These identities involve Euler numbers and are obtained by computing the moments of the random variable and the moments of the sum of two independent such random variables.
\end{abstract}

\maketitle

\markboth{\centerline{\scriptsize Probabilistic proof of a summation formula}}
{\centerline{\scriptsize T. Kim, D. S. Kim }}

\section{Introduction}

The aim of this paper is twofold. Firstly, we assume that $X$ is the hyperbolic secant random variable. Then we derive the summation formula in \eqref{7} for $\sum_{k=0}^{\infty} \frac{(-1)^{k}} {(2k+1)^{2n+1}}$ by computing the moments $E[X^{2n}]$ in two different ways. Secondly, we assume that $X$ and $Y$ are independent hyperbolic secant random variables. Then we obtain the expression in \eqref{40} for $\zeta(2n+2)$ by computing the moments $E[(X+Y)^{2n}]$ in two different ways. We note here that both of these identities involve the Euler numbers (see \eqref{1}). \par
In more detail, the outline of this paper is as follows. In Section 1, we first recall that Kim derived the explicit expression in \eqref{7} from the generating function of Euler numbers and the Fourier series of $f(x)=\sin\, ax$ on $[-\pi,\pi]$,\,(see [4]). Here $E_{n}$ are the Euler numbers given by \eqref{1}. We remind the reader of the previous result on the values of $\zeta(2n)$ involving $E_{2n-1}^{*}$,\, (see [4]). Here $E_{n}^{*}$ are the Euler numbers defined by $\frac{2}{e^t+1}=\sum_{n=0}^{\infty}E_{n}^{*}\frac{t^n}{n!}$. Let $X$ and $Y$ be independent continuous random variables, with $f(x)$, $g(x)$ their respective probability density functions. Then we recall that the probability density function of $X+Y$ is given by the convolution of $f(x)$ and $g(x)$. We remind the reader of the hyperbolic secant random variable. Finally, we recall the gamma function, its reflection formula and the beta function.
Section 2 contains the main results of this paper. Let $X$ be the hyperbolic secant random variable. We determine the moment generating function of $X$ by using beta function and the reflection formula of the gamma function. This yields that $E[X^{2n}]=E_{2n}(\frac{\pi}{2})^{2n}(-1)^{n}$. In Theorem 2, we obtain the same result as in \eqref{7} by computing the moment $E[X^{2n}]$ in another way directly from definition. Assume that $X$ and $Y$ are independent hyperbolic secant random variables. On the one hand, we derive the moment generating function of $X+Y$ from those of $X$ and $Y$, which were already determined. Thereby we obtain an expression of $E[(X+Y)^{2n}]$. On the other hand, we compute $E[(X+Y)^{2n}]$ directly from the definition. Now, equating these two give us an expression for $\zeta(2n+2),\,\,(n=0,1,2,\dots)$, in Theorem 4. In Section 3, we apply the central limit theorem to a sequence of independent hyperbolic secant random variables to show $\frac{\pi}{2} \sqrt{n} f_{X_1+X_2+\cdots+X_n} (\frac{\pi}{2} \sqrt{n}\, y) \rightarrow \frac{1}{\sqrt{2\pi}} e^{-\frac{y^2}{2}}$, \,\,as $n \rightarrow \infty$. Here $f_{X_1+X_2+\cdots+X_n}(x)$ is the probability density function of $X_1+X_2+\cdots+X_n$. In the rest of this section, we recall the facts that are needed throughout this paper.

\vspace{0.1in}

The Euler numbers are defined by
\begin{equation}
\sech x = \frac{1}{\cosh x} = \frac{2}{e^x +e^{-x}} =\sum_{n=0}^{\infty} E_n \frac{t^{n}}{n!},\quad |x| < \frac{\pi}{2}.\label{1}
\end{equation}
From \eqref{1}, we get
\begin{align}
& E_0 = 1,\, E_2 = -1,\, E_4 = 5,\, E_6 = -61,\, E_8=1385,\, E_{10}=-50521,\cdots;\nonumber \\
& E_{2k+1} = 0,\,\,\mathrm{for}~~ k = 0,1,2,\cdots, \nonumber  \quad (\mathrm{see}\ [4]).
\end{align} \par
Euler's formula states that, for any real number $x$,
\begin{equation}
e^{ix} = \cos x + i \sin x, \quad \mathrm{where} ~i =\sqrt{-1}, \quad(\mathrm{see}\ [13]). \label{2}
\end{equation}
From \eqref{2}, we note that
\begin{equation}
\cos x = \frac{1}{2} \big(e^{ix} + e^{-ix} \big), \quad (\mathrm{see}\ [1-11,13]). \label{3}	
\end{equation}
Thus, by \eqref{1} and \eqref{3}, we get
\begin{align}
\sec x & = \frac{2}{e^{ix} +e^{-ix}} = \sech (ix) = \sum_{n=0}^{\infty} \frac{i^n E_n}{n!} x^n \label{4}\\
& = \sum_{n=0}^{\infty} \frac{(-1)^n E_{2n}}{(2n)!} x^{2n} + i \sum_{n=0}^{\infty} \frac{(-1)^n E_{2n+1}}{(2n+1)!} x^{2n+1} \nonumber \\
& = \sum_{n=0}^{\infty} \frac{(-1)^n E_{2n}}{(2n)!} x^{2n}. \nonumber
\end{align}
From \eqref{4}, we have
\begin{equation}
x \sec x = \sum_{n=0}^{\infty} \frac{(-1)^n E_{2n}}{(2n)!} x^{2n+1},\quad \big(|x| < \frac{\pi}{2} \big).\label{5}
\end{equation} \par
From the Fourier series of $f(x) = \sin ax$ on $[-\pi, \pi ]$, Kim derived the following formula
\begin{equation}
\frac{\pi a}{2} \sec (\frac{\pi a}{2})  = \sum_{k=0}^{\infty}  2 \sum_{n=1}^{\infty} \frac{(-1)^{n-1}}{(2n-1)^{2k+1}} a^{2k+1},  \quad (\mathrm{see}\ [4]). \label{6}                                      \end{equation}
Thus, by \eqref{5} and \eqref{6}, we get
\begin{equation}
\sum_{k=0}^{\infty} \frac{(-1)^{k}} {(2k+1)^{2n+1}} = (-1)^n \frac{1}{2} \frac{E_{2n}} {(2n)!} \bigg(\frac{\pi}{2} \bigg)^{2n+1}=\frac{1}{2} \frac{|E_{2n}|} {(2n)!} \bigg(\frac{\pi}{2} \bigg)^{2n+1},\, (\mathrm{see}\ [4]), \label{7}
\end{equation}
where $n$ is a nonnegative integer. Here we note that Euler considered the numbers $|E_{2n}|$ in connection with sums like \eqref{7}. Later, in 1851, Raabe introduced the term ``Euler numbers,"\,\,(see [12]).

For $s \in \mathbb{C}$ with $Re(s) > 1$, the Riemann Zeta function is defined by
\begin{equation}
\zeta(s) = \sum_{n=1}^{\infty} \frac{1}{n^s},\quad (\mathrm{see}\ [1-11,13]). \label{8}
\end{equation}
From \eqref{7} and \eqref{8}, we note that
\begin{equation}
\zeta(2n+1,\frac{1}{4}) + 2^{2n} \big(1-2^{2n+1} \big) \zeta(2n+1) = (-1)^n \frac{E_{2n}}{2(2n)!}
\pi^{2n+1} 2^{2n}, \label{9}
\end{equation}
where $\zeta(s,a)$ is Hurwitz zeta function given by
\begin{equation*}
\zeta(s,a) = \sum_{n=0}^{\infty} \frac{1}{(n+a)^s}, \quad (a \neq 0,-1,-2, \cdots), \quad (\mathrm{see}\ [13]).
\end{equation*}
In addition, by \eqref{8}, we get
\begin{equation*}
\zeta(2n) =  \frac{(-1)^{n-1} (2\pi)^{2n} } {4(2n-1)! (1-4^n)}  E_{2n-1}^{*}
,\quad (\mathrm{see}\ [4]),
\end{equation*}
where $E_{n}^{*}$ is defined by $\frac{2}{e^t +1} = \sum_{n=0}^{\infty}
E_{n}^{*} \frac{t^n}{n!}.$ \par

Let $X$ and $Y$ be independent random variables such that $f(x)$ and $g(x)$ are their respective probability density functions. We recall that the cumulative distribution function of the random variable $X$ is defined by
\begin{equation*}
F_{X}(a) = P\{X \le a \} = \int_{-\infty}^{a} f(x) dx, \quad (\mathrm{see}\ [10]).
\end{equation*}
Assume that $f_{X+Y} (a)$ is the probability density function of
$X+Y$. Then it is given by the convolution of $f(x)$ and $g(x)$ as in the following:
\begin{equation}
f_{X+Y}(a) = \int_{-\infty}^{\infty} g(y) f(a-y) dy, \quad (\mathrm{see}\ [10]). \label{11}
\end{equation}

By \eqref{4}, we easily get
\begin{align}
\frac{1}{\pi} \int_{-\infty}^{\infty} \frac{1}{\cosh x} dx & = \frac{2}{\pi} \int_{0}^{\infty} \frac{1}{\cosh x} dx = \frac{4}{\pi} \int_{0}^{\infty} \frac{1}{e^x(1+e^{-2x})} dx \label{12}\\
& = \frac{4}{\pi} \sum_{k=0}^{\infty} (-1)^k \int_{0}^{\infty} e^{-(2k+1)x} dx
=  \frac{4}{\pi} \sum_{k=0}^{\infty} \frac{(-1)^k} {2k+1} = \frac{4}{\pi} \frac{\pi}{4} = 1. \quad \nonumber
\end{align}\par

A random variable $X$ is the hyperbolic secant random variable if the probability density function
is given by
\begin{equation}
f(x) = \frac{1}{\pi} \sech x = \frac{1}{\pi} \frac{1}{\cosh x}, \quad (x \in (-\infty, \infty)), \quad (\mathrm{see}\ [3]). \label{13}
\end{equation}

For $\alpha > 0$, the gamma function is defined by
\begin{equation}
\Gamma(\alpha) = \int_{0}^{\infty} e^{-t} t^{\alpha-1} dt, \quad (\mathrm{see}\ [13]). \label{14}
\end{equation}

For $\alpha, \beta > 0$, the Beta function is defined by
\begin{equation}
B(\alpha,\beta) = \int_{0}^{1} t^{\alpha-1} (1-t)^{\beta-1} dt, \quad (\mathrm{see}\ [13]). \label{15}
\end{equation}
Thus, by \eqref{15}, we get
\begin{equation}
B(\alpha,\beta) = \frac{\Gamma(\alpha) \Gamma(\beta)}{\Gamma(\alpha + \beta)}=\int_{0}^{\infty} \frac{t^{\beta-1}}{(1+t)^{\alpha+\beta}}dt, \quad (\mathrm{see}\ [13]). \label{16}
\end{equation}

For $|t| < 1$, the Euler's reflection formula of the gamma function is given by
\begin{equation}
\Gamma(t) \Gamma(1-t) = \frac{\pi}{\sin \pi t}, \quad (\mathrm{see}\ [13]). \label{17}
\end{equation}\par

\section{Probabilistic proof of a summation formula}

In this section, we assume that $X$ is the hyperbolic secant random variable. Then the probability density function of $X$ is given by
\begin{equation}
f(x) = \frac{1}{\pi} \sech x = \frac{1}{\pi} \frac{1}{\cosh x}, \quad (x \in (-\infty, \infty)). \label{18}
\end{equation}

First, we consider the moment generating function of $X$. \\
For $|t| < 1$, we have
\begin{align}
E[e^{Xt}] & = \int_{-\infty}^{\infty} e^{xt} f(x) dx = \frac{1}{\pi} \int_{-\infty}^{\infty} e^{xt} \frac{1}{\cosh x} dx \label{19}\\
& = \frac{2}{\pi} \int_{-\infty}^{\infty} e^{xt} \frac{1}{e^x + e^{-x}} dx
= \frac{2}{\pi} \int_{-\infty}^{\infty} e^{xt} \frac{e^{-x}}{e^{2x}+1} e^{2x} dx. \quad \nonumber
\end{align}
From \eqref{19}, by making change of the variable $y = e^{2x}$, and using \eqref{16} and \eqref{17}, we get
\begin{align}
E[e^{Xt}] & = \frac{1}{\pi} \int_{-\infty}^{\infty} e^{xt} \frac{e^{-x}}{e^{2x}+1} 2 e^{2x}dx = \frac{1}{\pi} \int_{0}^{\infty} \frac{y^{\frac{t}{2} + \frac{1}{2}-1}} {1+y} dy
 \label{20}\\
& = \frac{1}{\pi} B\bigg(\frac{t}{2} + \frac{1}{2}, 1-\frac{t}{2} - \frac{1}{2} \bigg)   = \frac{1}{\pi} \frac{\Gamma(\frac{t}{2} + \frac{1}{2}) \Gamma(1-\frac{t}{2} - \frac{1}{2})}{\Gamma(1)} \nonumber \\
&= \frac{1}{\pi} \frac{\pi}{\sin(\frac{t}{2} + \frac{1}{2}) \pi}= \frac{1}{\cos\frac{\pi t}{2}}.  \nonumber
\end{align}
From \eqref{4} and \eqref{20}, we note that
\begin{align}
\sum_{n=0}^{\infty} E[X^n] \frac{t^{n}}{n!} & = E[e^{Xt}]  = \frac{1}{\cos\frac{\pi t}{2}}
= \frac{1}{\cosh \frac{\pi i t}{2}}
\label{21}\\
& = \sum_{n=0}^{\infty} E_{2n}  \big(\frac{\pi}{2}\big)^{2n} (-1)^n \frac{t^{2n}}{(2n)!}. \quad \nonumber
\end{align}

Therefore, by \eqref{21}, we obtain the following theorem.

\begin{theorem}
For $n\ge 0$, we have
\begin{equation}
E[X^{2n}] = E_{2n} \big(\frac{\pi}{2}\big)^{2n} (-1)^n, \label{22}
\end{equation}
and
\begin{equation*}
E[X^{2n+1}] = 0.
\end{equation*}
\end{theorem}

We recall that the variance of $X$ is given by
\begin{equation}
\sigma^2 = E[\big(X - E[X] \big)^2] = E[X^2] - (E[X])^2. \label{23}
\end{equation}

From \eqref{22} and \eqref{23}, we note that
\begin{equation}
\sigma^2 = E[X^2] - (E[X])^2 = \big(\frac{\pi}{2}\big)^{2} E_2 (-1) = \big(\frac{\pi}{2}\big)^{2}, \quad \mu = E[X] = 0. \label{23-1}
\end{equation}

On the other hand, by \eqref{13}, we get
\begin{align}
(-1)^n \big(\frac{\pi}{2}\big)^{2n} E_{2n} & =  E[X^{2n}]  = \int_{-\infty}^{\infty} x^{2n} f(x) dx \label{24}\\
& = \frac{1}{\pi} \int_{-\infty}^{\infty} x^{2n} \frac{1}{\cosh x} dx
= \frac{4}{\pi} \int_{0}^{\infty} x^{2n}\frac{1}{e^x + e^{-x}} dx  \nonumber \\
& = \frac{4}{\pi} \int_{0}^{\infty} x^{2n}\frac{e^{-x}}{1 + e^{-2x}} dx
  = \frac{4}{\pi} \sum_{k=0}^{\infty} (-1)^k \int_{0}^{\infty} x^{2n} e^{-(2k+1)x} dx. \quad \nonumber
\end{align}

Now, we observe that
\begin{align}
\int_{0}^{\infty} x^{2n} e^{-(2k+1)x} dx & =  \bigg[-\frac{e^{-(2k+1)x}}{2k+1} x^{2n} \bigg]_{0}^{\infty} + \frac{2n}{2k+1} \int_{0}^{\infty} x^{2n-1} e^{-(2k+1)x} dx  \label{25} \\
& = \frac{2n}{2k+1} \int_{0}^{\infty} x^{2n-1} e^{-(2k+1)x} dx
  = \frac{2n(2n-1)}{(2k+1)^2} \int_{0}^{\infty} x^{2n-2} e^{-(2k+1)x} dx \nonumber \\
& = \cdots \nonumber \\
& = \frac{2n(2n-1) \cdots 1}{(2k+1)^{2n}} \int_{0}^{\infty} e^{-(2k+1)x} dx = \frac{(2n)!}{(2k+1)^{2n+1}}.\nonumber
\end{align}
Therefore, by \eqref{24} and \eqref{25}, we get
\begin{equation}
(-1)^n \big(\frac{\pi}{2}\big)^{2n} E_{2n} = \frac{4}{\pi} \sum_{k=0}^{\infty} \frac{(-1)^k}{(2k+1)^{2n+1}} (2n)!. \label{26}
\end{equation}
From \eqref{26}, we have the following theorem.
\begin{theorem}
For $n\ge 0$, we have
\begin{displaymath}
\sum_{k=0}^{\infty} \frac{(-1)^k}{(2k+1)^{2n+1}} = \frac{(-1)^n}{(2n)!} \frac{1}{2} \big(\frac{\pi}{2}\big)^{2n+1} E_{2n}.
\end{displaymath}
\end{theorem}
In particular, for $n=1$, we have the following result.
\begin{corollary}
\begin{displaymath}
\sum_{k=0}^{\infty} \frac{(-1)^k}{(2k+1)^{3}} = -\frac{1}{4} \big(\frac{\pi}{2}\big)^{3} E_{2}
= \frac{1}{32} \pi^3.
\end{displaymath}
\end{corollary} \par

Assume that $X$ and $Y$ are independent hyperbolic secant random variables. \\
Then, by \eqref{11}, we get
\begin{align}
f_{X+Y}(a) & =  \int_{-\infty}^{\infty} \frac{1}{\pi} \frac{1}{\cosh(y)} \frac{1}{\pi} \frac{1}{\cosh(a-y)} dy \label{28}\\
& = \frac{4}{\pi^2} \int_{-\infty}^{\infty} \frac{1}{e^y + e^{-y}} \frac{1}{e^{a-y} + e^{y-a}} dy
 \nonumber \\
& = \frac{4}{\pi^2} \int_{-\infty}^{\infty} \frac{e^y}{e^{2y}+1} \frac{e^{-a+y}}{e^{-2a+2y} + 1} dy. \quad \nonumber
\end{align}

From \eqref{28}, by making change of the variable $e^{y} = x$, we have
\begin{align}
f_{X+Y}(a) & = \frac{4}{\pi^2} \int_{0}^{\infty} \frac{x}{1+x^2} \frac{e^{-a}}{1+x^2e^{-2a}} dx \label{29}\\
& = \frac{4}{\pi^2(e^a-e^{-a})} \int_{0}^{\infty} \bigg(\frac{x}{1+x^2} - \frac{xe^{-2a}}{1+x^2e^{-2a}} dx  \nonumber \\
& = \frac{4}{\pi^2(e^a-e^{-a})} \bigg[\frac{\log(1+x^2) - \log(1+x^2)e^{-2a}}{2} \bigg]_{0}^{\infty} \nonumber \\
& = \frac{4}{\pi^2(e^a-e^{-a})} \bigg[\frac{1}{2} \log\big(\frac{1+x^2}{1+x^2e^{-2a}}\big) \bigg]_{0}^{\infty}\nonumber \\
& = \frac{4}{\pi^2(e^a-e^{-a})} \frac{1}{2} \log(e^{2a}) = \frac{4a}{\pi^2(e^a-e^{-a})}, \quad \nonumber
\end{align}
where $a \in (-\infty, \infty)$.

Thus, the probability density function of $Z = X+Y$ is given by
\begin{equation}
f_Z (a) = f_{X+Y} (a) = \frac{4a}{\pi^2 (e^a - e^{-a})}, \quad (a \in (-\infty, \infty)). \label{30}
\end{equation}
Since $X$ and $Y$ are independent random variables,
\begin{equation}
E[e^{t(X+Y)}]  = E[e^{tX}] E[e^{tY}] = \frac{1}{\cos \frac{\pi}{2} t} \frac{1}{\cos \frac{\pi}{2} t}. \label{31}
\end{equation}
Thus, by \eqref{4} and \eqref{31}, we get
\begin{align}
E[e^{t(X+Y)}] & = \sum_{m=0}^{\infty} E_{2m}  \big(\frac{\pi}{2}\big)^{2m} (-1)^m \frac{t^{2m}}{(2m)!} \sum_{j=0}^{\infty} E_{2j}  \big(\frac{\pi}{2}\big)^{2j} (-1)^j \frac{t^{2j}}{(2j)!} \label{32}\\
& = \sum_{n=0}^{\infty} \bigg( (-1)^n \big(\frac{\pi}{2}\big)^{2n} \sum_{l=0}^{n} \binom{2n}{2l} E_{2l}
E_{2n-2l} \bigg) \frac{t^{2n}}{(2n)!}, \quad \nonumber
\end{align}
From \eqref{32}, we note that
\begin{align}
(-1)^n \big(\frac{\pi}{2}\big)^{2n} \sum_{l=0}^{n} \binom{2n}{2l} E_{2l}
E_{2n-2l} & = E[(X+Y)^{2n}] = E[Z^{2n}] \label{33}\\
& = \int_{-\infty}^{\infty} x^{2n} f_Z (x) dx = \int_{-\infty}^{\infty} x^{2n} \frac{4x}{\pi^2} \frac{1}{e^x - e^{-x}} dx \nonumber \\
& = \frac{8}{\pi^2} \int_{0}^{\infty} \frac{x^{2n+1}}{e^x - e^{-x}} dx = \frac{8}{\pi^2} \int_{0}^{\infty} \frac{e^{-x} x^{2n+1}}{1 - e^{-2x}} dx \nonumber \\
& = \frac{8}{\pi^2} \sum_{k=0}^{\infty} \int_{0}^{\infty} x^{2n+1} e^{-(2k+1)x} dx. \quad \nonumber
\end{align}
Now, we observe that
\begin{align}
\int_{0}^{\infty} x^{2n+1} e^{-(2k+1)x} dx & =  \bigg[-\frac{x^{2n+1}}{2k+1}  e^{-(2k+1)x} \bigg]_{0}^{\infty} + \frac{2n+1}{2k+1} \int_{0}^{\infty} x^{2n} e^{-(2k+1)x} dx  \label{34} \\
& = \frac{2n+1}{2k+1} \int_{0}^{\infty} x^{2n} e^{-(2k+1)x} dx= \frac{(2n+1) 2n}{(2k+1)^2} \int_{0}^{\infty} x^{2n-1} e^{-(2k+1)x} dx \nonumber \\
& = \cdots \nonumber \\
& = \frac{(2n+1) 2n \cdots 2 \cdot 1}{(2k+1)^{2n+1}} \int_{0}^{\infty} e^{-(2k+1)x} dx = \frac{(2n+1)!}{(2k+1)^{2n+2}}.\nonumber
\end{align}
By \eqref{33} and \eqref{34}, we get
\begin{equation}
(-1)^n \big(\frac{\pi}{2}\big)^{2n} \sum_{l=0}^{n} \binom{2n}{2l} E_{2l}
E_{2n-2l} = \frac{8}{\pi^2} (2n+1)! \sum_{k=0}^{\infty} \frac{1}{(2k+1)^{2n+2}}. \label{35}
\end{equation}
By \eqref{35}, we have
\begin{align}
\sum_{k=0}^{\infty} \frac{1}{(2k+1)^{2n+2}} & =  (-1)^n \big(\frac{\pi}{2}\big)^{2n} \frac{\pi^2}{8}
\frac{1}{(2n+1)!} \sum_{l=0}^{n} \binom{2n}{2l} E_{2l}
E_{2n-2l} \label{36} \\
& = (-1)^n \big(\frac{\pi}{2}\big)^{2n+2} \frac{1}{2}
\frac{1}{(2n+1)!} \sum_{l=0}^{n} \binom{2n}{2l} E_{2l}
E_{2n-2l}.\nonumber
\end{align}
By \eqref{8}, we have
\begin{align}
\zeta(2n+2) & =  \sum_{k=1}^{\infty} \frac{1}{k^{2n+2}} \label{37} \\
& =  \sum_{k=0}^{\infty} \frac{1}{(2k+1)^{2n+2}} + \sum_{k=1}^{\infty} \frac{1}{(2k)^{2n+2}} \nonumber \\
& =  \sum_{k=0}^{\infty} \frac{1}{(2k+1)^{2n+2}} + \frac{1}{2^{2n+2}} \zeta(2n+2).\nonumber
\end{align}
Thus, \eqref{37}, we get
\begin{equation}
\sum_{k=0}^{\infty} \frac{1}{(2k+1)^{2n+2}} = \bigg( 1 - \big(\frac{1}{2}\big)^{2n+2} \bigg) \zeta(2n+2). \label{38}
\end{equation}
From \eqref{36} and \eqref{38}, we have the following theorem.
\begin{theorem}
	For $n\ge 0$, we have
	\begin{equation}
	\zeta(2n+2) = \frac{(-1)^n} {\bigg( 1 - \big(\frac{1}{2}\big)^{2n+2} \bigg) (2n+1)!} \frac{1}{2} \big(\frac{\pi}{2}\big)^{2n+2} \sum_{l=0}^{n} \binom{2n}{2l} E_{2l}
	E_{2n-2l}. \label{40}
	\end{equation}	
\end{theorem}
In particular, if $n=0$ in \eqref{40}, then we get
\begin{displaymath}
\zeta(2) = \frac{4}{3} \frac{1}{2} \big(\frac{\pi}{2}\big)^{2} = \frac{\pi^2}{6}.
\end{displaymath}

\section{Further Remark}
We first recall the following Central Limit Theorem (see [10]).
\begin{lemma} [Central Limit Theorem]
Let $X_1, X_2, \dots$\, be a sequence of independent random variables with $E[X_i] = \mu, ~ Var(X_i) = \sigma^2,\quad (i = 1,2, \dots).$\, Then we have
\begin{displaymath}
\lim_{n \rightarrow \infty} \frac{X_1+X_2+\cdots+X_n - n\mu}{\sigma \sqrt{n}} \sim N(0,1).
\end{displaymath}
That is,
\begin{displaymath}
P\bigg\{\frac{X_1+X_2+\cdots+X_n - n\mu}{\sigma \sqrt{n}}
\le a \bigg\} \rightarrow \frac{1}{\sqrt{2\pi}} \int_{-\infty}^{a} e^{-\frac{x^2}{2}} dx,
\end{displaymath}
 as $n \rightarrow \infty.$
\end{lemma}

Let $X_1, X_2, \dots$ be a sequence of independent hyperbolic secant random variables.
Then by \eqref{23-1}, we get
\begin{equation}
\sigma^2 = Var(X_i) = \big(\frac{\pi}{2}\big)^2, ~  \mu = E[X_i] = 0, \quad  \label{41}
\end{equation}
where $i = 1,2, \dots.$

By Central Limit Theorem, we have
\begin{equation}
P\bigg\{\frac{X_1+X_2+\cdots+X_n}{\frac{\pi}{2} \sqrt{n}}
\le y \bigg\} \rightarrow \frac{1}{\sqrt{2\pi}} \int_{-\infty}^{y} e^{-\frac{x^2}{2}} dx,\label{42}
\end{equation}
as $n \rightarrow \infty.$

Let $f_{X_1+X_2+\cdots+X_n} (x)$ be the probability density function of $X_1+X_2+\cdots+X_n$.
Then \eqref{42} is the same as saying that
\begin{equation}
\int_{-\infty}^{y}\frac{\pi}{2}\sqrt{n}f_{X_1+X_2+\cdots+X_n}(\frac{\pi}{2}\sqrt{n} \,x) dx \rightarrow \int_{-\infty}^{y} \frac{1}{\sqrt{2 \pi}}e^{-\frac{x^2}{2}}dx, \label{43}
\end{equation}
as $n \rightarrow \infty$.
Thus, by \eqref{43}, we get
\begin{displaymath}
\frac{\pi}{2} \sqrt{n} f_{X_1+X_2+\cdots+X_n} (\frac{\pi}{2} \sqrt{n}\, y) \rightarrow \frac{1}{\sqrt{2\pi}} e^{-\frac{y^2}{2}},
\end{displaymath}
 as $n \rightarrow \infty.$

\section{Conclusion}
In this paper, among other things, we showed the summation formula
\begin{equation*}
\sum_{k=0}^{\infty} \frac{(-1)^{k}} {(2k+1)^{2n+1}} = (-1)^n \frac{1}{2} \frac{E_{2n}} {(2n)!} \bigg(\frac{\pi}{2} \bigg)^{2n+1}=\frac{1}{2} \frac{|E_{2n}|} {(2n)!} \bigg(\frac{\pi}{2} \bigg)^{2n+1},
\end{equation*}
by evaluating even moments of the hyperbolic secant random variable in two different ways.
As we mentioned in the Introduction, this formula was known to Euler. So the contribution of the present paper would be that it gives a probabilistic and simple proof of the above summation formula. \par
We would like to note that the Euler numbers in \eqref{1} has interesting connections with Euler zigzag numbers defined by the Taylor series, which is given by
\begin{equation*}
\sec x+\tan x=\tan \Big(\frac{\pi}{4}+\frac{x}{2} \Big)=\sum_{n=0}^{\infty}A_n \frac{x^n}{n!}.
\end{equation*}
Indeed, one can show that $A_{2n}=(-1)^n E_{2n}=|E_{2n}|,\,\, (n=0,1,2,\dots)$. Interesting combinatorial interpretations for the zigzag numbers can be found in the stanley's  recent talk given at the 14th Ramanujan Colloquium, which was held at University of Florida in 2023,\, (see [12]). For example, $A_n$ is equal to the number of permutations in the symmetric group $S_{n}$ that are alternationg. Here a sequence $a_1,a_2,\dots,a_n$ of distinct integers is defined to be {\it{alternating}} if $a_1 > a_2 < a_3 > a_4 < \cdots$. For example, $A_4=E_4=5$, since $2143,\,3142,\,3241,\,4132,\,4231$ are the alternating permutations in $S_4$.

\end{document}